
%

\baselineskip=14pt
\parskip=10pt

\magnification=\magstephalf

\def\1{{\overline{1}}}
\def\2{{\overline{2}}}
\parindent=0pt
\overfullrule=0in

\def\frac#1#2{{#1 \over #2}}

\bf
\centerline
{
The  Generating Functions Enumerating ${\bf 12...d}$-Avoiding Words with ${\bf r}$ occurrences  of each
}
\centerline
{
 of ${\bf 1 \, , \, \dots \, , \,  n }$ are  ${\bf D}$-finite for all ${\bf d}$ and all ${\bf r}$
}

\rm
\bigskip
\centerline
{\it By Shalosh B. EKHAD and Doron ZEILBERGER}
\bigskip
\qquad \qquad 
{\it 
Dedicated to Neil James Alexander Sloane (born October 10, 1939), on his \hfill\break A000027[75]-th birthday
(alias  A005408[38]-th, alias  A002808[53]-th, alias A001477[74]-th, alias A014612[17]-th,
and 13520 other aliases.)
}

\bigskip

{\bf Introduction}

In a recent beautiful article, Nathaniel Shar and Doron Zeilberger ([ShZ]) 
proved that for any positive integer $r$,
the generating function of the sequence enumerating $123$-avoiding words with $r$ occurrences of each of the 
letters $1, \, \dots \, , n$ is always {\bf algebraic}.  In other words for each $r$, the generating
function, let's call it $f_r(x)$, satisfies an equation of the form
$$
P_r(x, f_r(x))=0 \quad,
$$
for {\it some} polynomial, $P_r$, of two variables. The actual polynomials, $P_r(x,y)$, were computed for $r \leq 4$.

This is no longer true for $12 \dots d$-avoiding words with $d\geq 4$, {\it even} for $r=1$.

In 1990 Doron Zeilberger ([Z]) showed that for and {\it any} positive integer $d$, the generating function
enumerating $1\dots d$-avoiding {\it permutations} (i.e. words in $\{1, \, \dots \, , n \}$ where each letter occurs exactly $1$ times)
is the next-best-thing to being algebraic, which is being  {\bf D-finite} (aka as {\it holonomic}).
Recall that a formal power series is $D$-finite if it satisfies a linear {\it differential} equation
with polynomial coefficients, or equivalently, the enumerating sequence itself is {\bf P-recursive}, i.e. 
satisfies a linear {\it recurrence} equation with polynomial coefficients. 
Ira Gessel ([G]) famously discovered (and proved) a beautiful determinant with Bessel functions, 
for the generating function, (of the sequence divided by $n!^2$) (that also implies the above result), 
and Amitai Regev ([R]) famously derived delicate and precise asymptotics.

In the present article, dedicated to guru Neil Sloane on his 75-th birthday, we observe that the analogous 
generating functions for {\it multi-set permutations} (alias {\it words}), where every letter appears the same number of times, say $r$, 
are still always {\bf D-finite}, (for {\it every} $d$ and {\it every} $r$), 
and we actually crank out the first few terms of quite a few of them, 
many of whom are not yet in the OEIS ([Sl]). 

All this data, often with linear recurrences (that we know {\it a priori} exist, and hence
it justifies their discovery by pure guessing), and very precise asymptotics, is collected in the
front of this article

{\tt http://www.math.rutgers.edu/\~{}zeilberg/mamarim/mamarimhtml/sloane75.html}

where links to two useful Maple packages, that were used to generate all that data,
{\tt SLOANE75} and {\tt NEIL}, can be found and downloaded, and readers who have Maple
and computer time to spare are welcome to use in order to generate yet more data.

Last but not least, we pledge  $100$ dollars to the OEIS in honor of the first one to 
prove our conjectured asymptotic formula for the number of $1 \dots d$-avoiding words in $\{1^r \dots n^r\}$
that generalizes Regev's ([R]) famous formula for $r=1$.
We pledge  another $100$ dollars for extending Gessel's Bessel determinant, from the $r=1$ case to general $r$.

{\bf Why is the Sequence Enumerating ${\bf 1 \dots d}$-avoiding words in ${\bf \{1^r \dots n^r\}}$ P-recursive?}

By the Robinson-Schenstead-Knuth (RSK) famous correspondence, our quantity of interest,
let's call it $A_{d,r}(n)$ is given by
$$
A_{d,r}(n) \, = \,
\sum_{ {\lambda \vdash rn} \atop  {length(\lambda) \leq d} } \, f_{\lambda} g^{(r)}_{\lambda} \quad,
$$
where $f_{\lambda}$ is the number of standard Young tableaux of shape $\lambda=(\lambda_1, \dots, \lambda_d)$, and
$g^{(r)}_{\lambda}$ is the number of column-strict Young tableaux with exactly $r$ occurrences of each of $1,\, \dots \, , n$.
For $\lambda=(\lambda_1, \dots, \lambda_d)$ (where we pad it with zeroes if the length is less than $d$),
$f_{\lambda}$ is closed-form (thanks to Young-Frobenius, or the hook-length formula), and hence {\it ipso facto},
holonomic in its $d$ discrete arguments. Furthermore, for $r>1$, 
$g^{(r)}_{\lambda}$, while no longer closed-form, is easily seen to be holonomic in its $d$ discrete arguments
(one way to see this is to note  that their {\it redundant} generating function
(in the sense of MacMahon) is a rational formal power series in $x_1, \dots, x_d$).
It follows, by {\it general holonomic nonsense} ([Z]), that for any {\bf fixed} integers $r$ and $d$
the sequence, in $n$, $\{ A_{d,r}(n) \}$, is $P$-recursive. Computationally speaking, 
it is fairly easy to compute  $g^{(r)}_{\lambda}$, and hence crank-out the first few terms of
the sequences  $\{ A_{d,r}(n) \}$ for quite a few $d$ and $r$, that for $r$ and $d$ not too large
may be used to guess (in {\it real} time) the recurrences empirically, that we know
must be the right ones.

{\bf The 100 dollars conjecture generalizing Regev's Asymptotics}

{\bf Conjecture} ($100$ donation to the OEIS in honor of the first prover)

Let $A_{d,r}(n)$ be the number of $1 \dots d$-avoiding words in $\{1^r \dots n^r\}$, then
there exists a constant $C_{r,d}$ such that
$$
A_{d,r}(n) \, \sim \, C_{r,d} \cdot \left ( { { d+r-2} \choose {d-2} } (d-1)^r \right )^n \cdot \frac{1}{n^{((d-1)^2-1)/2}} \quad .
$$

{\bf Extra Credit} ($25$ additional dollars): find an explicit expression for $C_{r,d}$ in terms of $r$ and $d$ (involving $\pi$, of course).

{\bf The 100 dollars Challenge to generalize Gessel's Spectacular Theorem}

This is more open-ended, but it would be nice to get a determinant expression, in the style of Ira Gessel's ([G])
famous expression for the generating function of $A_{d,1}(n)/n!^2$, {\it canonized} in the {\it bible} ([W], p. 996, Eq. (5)).
Here it is: Let $u_k(n):=A_{k,1}(n)$, then
$$
\sum_{n \geq 0} \frac{u_k(n)}{n!^2} \, x^{2n} = \, \det (I_{|i-j|}(2x))_{i,j=1, \dots, k}   \quad ,
$$
in which $I_{\nu}(t)$ is (the modified Bessel function)
$$
I_\nu (t) = \sum_{j=0}^{\infty} \frac{ (\frac{1}{2} \, t)^{2j+\nu}} {j!(j+\nu)!} \quad .
$$

Guru Herb Wilf (ibid) goes on to wax eloquently:

\quad \quad 
{\it
`` At any rate, it seems fairly ``spectacular'' to me that when you place various infinite series
such as the above into a $k \times k$ determinant, and then expand the determinant, you should find that
the coefficient of $x^{2n}$, when multiplied by $n!^2$, is exactly the number of permutations
of $n$ letters with no increasing subsequence longer than $k$ .''}

It would be even more spectacular, if you, dear reader, would generalize this to $r>1$!

{\bf References}

[G] I.  Gessel, {\it Symmetric functions and P-recursiveness},
Journal of Combinatorial Theory, Series A {\bf 53} (1990), 257-285; \hfill\break
{\tt http://people.brandeis.edu/\~{}gessel/homepage/papers/dfin.pdf} \quad .

[R] A. Regev, {\it Asymptotic values for degrees associated with strips of Young diagrams},
Adv. Math. {\bf 41} (1981), 115-136.

[ShZ]  N. Shar and D. Zeilberger,
{\it The (ordinary) generating functions enumerating $123$-avoiding words with $r$ occurrences of each of $1,2, ..., n$ are always algebraic},
submitted; \hfill\break
{\tt http://www.math.rutgers.edu/\~{}zeilberg/mamarim/mamarimhtml/words123.html} \quad .

[Sl] N. J.A. Sloane, The Online-Encyclopedia of Integer Sequences (OEIS) {\it https://oeis.org/} \quad .

[W] H. Wilf,  {\it Mathematics, an experimental science}, in: 
``Princeton Companion to Mathematics'', (W. Timothy Gowers, ed.), Princeton University Press, 2008, 991-1000; \hfill\break
{\tt http://www.math.rutgers.edu/\~{}zeilberg/akherim/HerbMasterpieceEM.pdf} \quad .

[Z]  D. Zeilberger, {\it A Holonomic Systems approach to Special Functions}, J. Computational and Applied Math {\bf 32} (1990), 321-368;
\hfill\break
{\tt http://www.math.rutgers.edu/\~{}zeilberg/mamarim/mamarimhtml/holonomic.html} \quad .

\vfill\eject

\hrule
\bigskip
Doron Zeilberger, Department of Mathematics, Rutgers University (New Brunswick), Hill Center-Busch Campus, 110 Frelinghuysen
Rd., Piscataway, NJ 08854-8019, USA. \hfill \break
zeilberg at math dot rutgers edu \quad ;  \quad {\tt http://www.math.rutgers.edu/\~{}zeilberg/} \quad .
\bigskip
\hrule
\bigskip
Shalosh B. Ekhad, c/o D. Zeilberger, Department of Mathematics, Rutgers University (New Brunswick), Hill Center-Busch Campus, 110 Frelinghuysen
Rd., Piscataway, NJ 08854-8019, USA.
\bigskip
\hrule

\bigskip
EXCLUSIVELY PUBLISHED IN  The Personal Journal of Shalosh B. Ekhad and Doron Zeilberger ({ \tt http://www.math.rutgers.edu/\~{}zeilberg/pj.html})
 and {\tt arxiv.org}.
\bigskip
\hrule
\bigskip
{\bf Dec. 5, 2014}

\end